\documentclass[a4paper,11pt]{article}

\usepackage[matrix,arrow,curve]{xy}

\usepackage{latexsym}
\usepackage{amsmath}
\usepackage{amssymb}
\usepackage[mathscr]{eucal}

\newtheorem{theorem}{Theorem}[section]
\newtheorem{proposition}[theorem]{Proposition}
\newtheorem{lemma}[theorem]{Lemma}
\newtheorem{corollary}[theorem]{Corollary}
\newtheorem{definition}[theorem]{Definition}
\newtheorem{example}[theorem]{Example}
\newtheorem{remark}[theorem]{Remark}

\newcommand{\brem}{\begin{remark}\rm}
\newcommand{\erem}{\end{remark}}
\newcommand{\bex}{\begin{example}\rm}
\newcommand{\eex}{\end{example}}
\newcommand{\bthm}{\begin{theorem}}
\newcommand{\ethm}{\end{theorem}}
\newcommand{\blem}{\begin{lemma}}
\newcommand{\elem}{\end{lemma}}
\newcommand{\bprop}{\begin{proposition}}
\newcommand{\eprop}{\end{proposition}}
\newcommand{\bcor}{\begin{corollary}}
\newcommand{\ecor}{\end{corollary}}
\newcommand{\bdefi}{\begin{definition}}
\newcommand{\edefi}{\end{definition}}

    \newcommand{\qedsymbol}{~\rule[-0.35mm]{2mm}{2mm}}
\newenvironment{env}[2]{\begin{#1}#2\end{#1}}{}
    \newcommand{\beq}[1]{\begin{env}{equation}{#1}}
    \newcommand{\beqn}[1]{\begin{env}{equation*}{#1}}
    \newcommand{\bal}[1]{\begin{env}{align}{#1}}
    \newcommand{\baln}[1]{\begin{env}{align*}{#1}}
    \newcommand{\bga}[1]{\begin{env}{gather}{#1}}
    \newcommand{\bgan}[1]{\begin{env}{gather*}{#1}}
    \newcommand{\bflal}[1]{\begin{env}{flalign}{#1}}
    \newcommand{\bflaln}[1]{\begin{env}{flalign*}{#1}}
    \newcommand{\bmu}[1]{\begin{env}{multline}{#1}}
    \newcommand{\bmun}[1]{\begin{env}{multline*}{#1}}
    \newcommand{\bsp}[1]{\begin{env}{split}{#1}}

    \newcommand{\eeq}{\end{env}}
    \newcommand{\eeqn}{\end{env}}
    \newcommand{\eal}{\end{env}}
    \newcommand{\ealn}{\end{env}}
    \newcommand{\ega}{\end{env}}
    \newcommand{\egan}{\end{env}}
    \newcommand{\eflal}{\end{env}}
    \newcommand{\eflaln}{\end{env}}
    \newcommand{\emu}{\end{env}}
    \newcommand{\emun}{\end{env}}
    \newcommand{\esp}{\end{env}}

\newcommand{\lf}{\vspace{2ex}}

\renewcommand{\bf}[1]{\textbf{#1}}
\renewcommand{\it}[1]{\textit{#1}}
\renewcommand{\sc}[1]{\textsc{#1}}
\renewcommand{\sf}[1]{\textsf{#1}}

\renewcommand{\tt}[1]{\texttt{#1}}
\newcommand{\hl}[1]{\bf{\it{#1}}}
\newcommand{\mrm}[1]{\mathrm{#1}}
\newcommand{\mbf}[1]{\mathbf{#1}}

\newcommand{\cmc}[1]{\mathcal{#1}}
\newcommand{\eus}[1]{\mathscr{#1}}

\newcommand{\bb}[1]{\mathbb{#1}}
\newcommand{\nbd}[1]{$#1$\nobreakdash--}
\newcommand{\ol}[1]{\overline{#1}}

\newcommand{\vt}{\vartheta}

\newcommand{\vp}{\varphi}

\newcommand{\abs}[1]{\left\lvert#1\right\rvert}

\newcommand{\bfam}[1]{\bigl(#1\bigr)}
\newcommand{\Bfam}[1]{\Bigl(#1\Bigr)}
\newcommand{\AB}[1]{\langle#1\rangle}
\newcommand{\bAB}[1]{\bigl\langle#1\bigr\rangle}

\newcommand{\CB}[1]{\{#1\}}
\newcommand{\bCB}[1]{\bigl\{#1\bigr\}}
\newcommand{\BCB}[1]{\Bigl\{#1\Bigr\}}
\newcommand{\SB}[1]{[#1]}

\newcommand{\set}[2][]{
    \ifthenelse{\equal{#1}{}}{
        \CB{#2}}{
        \CB{#1~|~#2}}}
\newcommand{\bset}[2][]{
    \ifthenelse{\equal{#1}{}}{
        \bCB{#2}}{
        \bCB{#1~|~#2}}}
\newcommand{\Bset}[2][]{
    \ifthenelse{\equal{#1}{}}{
        \BCB{#2}}{
        \BCB{#1~\big|~#2}}}
\newcommand{\zero}{\CB{0}}

\DeclareMathOperator{\ls}{\normalfont\sf{span}}
\DeclareMathOperator{\cls}{\ol{\ls}}

\DeclareMathOperator{\id}{\normalfont\sf{id}}

\newcommand{\C}{\bb{C}}

\newcommand{\E}{\bb{E}}

\newcommand{\N}{\bb{N}}

\newcommand{\cA}{\cmc{A}}
\newcommand{\cB}{\cmc{B}}
\newcommand{\cC}{\cmc{C}}
\newcommand{\cD}{\cmc{D}}

\newcommand{\cK}{\cmc{K}}

\newcommand{\sB}{\eus{B}}

\newcommand{\U}{\mbf{1}}

\newcommand{\eins}{\mbox{$1\!\!\!\:{\rm I}$}}

\newcommand{\BG}{\mbox{$\sB(G)$}}

\newcommand{\BK}{\mbox{$\sB(K)$}}

    \numberwithin{equation}{section}

    \newcommand{\proofname}{Proof}
    \renewcommand{\qedsymbol}{~\rule[-0.35mm]{2mm}{2mm}}
    \newenvironment{Proof}{
        \vspace{1ex}}{
        \nolinebreak\qedsymbol}
    \newcommand{\proof}{
        \begin{Proof}{\proofname}\ignorespaces}
    \newcommand{\qed}{\end{Proof}}
    \newcommand{\noqed}{
        \renewcommand{\qedsymbol}{}
        \end{Proof}}

\begin{document}

\title{Constructing Extensions of CP-Maps\\via Tensor Dilations\\with the Help of Von Neumann Modules}
\author{Rolf Gohm \and Michael Skeide
\thanks{RG is supported by DFG, MS is supported by DAAD, ISI Bangalore and Research Fonds of the Department S.E.G.e S.\ of University of Molise.}}
\date{Campobasso and Greifswald, 4 November, 2003}
\maketitle

\begin{abstract}
\noindent
We apply Hilbert module methods to show that normal completely positive
maps admit weak tensor dilations. Appealing to a duality between weak
tensor dilations and extensions of CP-maps, we get an existence proof for
certain extensions. We point out that this duality is part of a far reaching
duality between a von Neumann bimodule and its commutant in which also other
dualities, known and new, find their natural common place.
\end{abstract}

\section{Introduction}

A \hl{dilation} of a CP-map (a completely positive  mapping) $S$ between $C^*$--al\-ge\-bras $\cA$ and $\cB$ is a homomorphism $j$ from $\cA$ to another $C^*$--algebra $\cD$ containing $\cB$ and a conditional expectation $P\colon\cD\rightarrow\cB$ such that $S$ is recovered as $S=P\circ j$.
\beqn{
\xymatrix{
\cA	\ar[drr]_j	\ar[rr]^S	&&	\cB						\\
					&&	**[r]	\cD\supset\cB	\ar[u]_P
}
}\eeqn
We assume that the $C^*$--algebras have a unit. A dilation is \hl{unital}, if $j$ is unital and if the unit of $\cB\subset\cD$ is the unit of $\cD$. On the contrary, a dilation is a \hl{weak} dilation, if the conditional expectation has the form $P(d)={\eins}_{\cB}\, d\, {\eins}_{\cB}$ (i.e.\ $\cB$ is a corner of $\cD$).
\beqn{
\xymatrix{
\cA	\ar[drr]_j	\ar[rr]^S	&&	\cB								\\
					&&	\cD	\ar[u]_{{\eins}_{\cB}\bullet{\eins}_{\cB}}
}
}\eeqn
This excludes (except in trivial cases) that ${\eins}_{\cB}={\eins}_{\cD}$.

\brem
The name weak dilation is borrowed from Bhat and Par\-tha\-sa\-ra\-thy 
\cite{BhPa94}. They consider the case when $\cA=\cB$ and when $j$ factors 
through a (usually unital) endomorphism $\vartheta$ of $\cD$ 
(i.e.\ $j=\vartheta\circ\id_{\cB}$) such that a simultaneous dilation is obtained 
for the whole semigroup $(S^n)_{n\in\N_0}$ to the semigroup 
$(\vartheta^n)_{n\in\N_0}$ of (unital) endomorphisms of $\cD$.
\beqn{
\parbox{3cm}{
\xymatrix{
\cB	\ar[d]_\subset	\ar[drr]_j	\ar[rr]^S	&&	\cB								\\
\cD	\ar[rr]_\vt							&&	\cD	\ar[u]_{{\eins}_{\cB}\bullet{\eins}_{\cB}}
}
}
~~~~~~\Longrightarrow~~~~~~
\parbox{3cm}{
\xymatrix{
\cB	\ar[d]_\subset	\ar[rr]^{S^n}	&&	\cB								\\
\cD	\ar[rr]_{\vt^n}					&&	\cD	\ar[u]_{{\eins}_{\cB}\bullet{\eins}_{\cB}}
}
}
}\eeqn
Such a dilation is then a weak dilation in the sense of \cite{BhPa94} and, 
in particular, it is a weak dilation in the sense above.\erem

Gohm \cite{Goh04,Goh03a} worked with \it{weak tensor dilations}, where the ``big'' algebra $\cD$ is just a tensor product $\cB\otimes\cC$, and where the conditional expectation is obtained from a state on the second factor. (We give immediately the version for von Neumann algebras, because our results apply only to that case. The tensor product is, therefore, that of von Neumann algebras.)

\begin{definition}
Let $\cA, \cB, \cC $ be von Neumann algebras. Let $S\colon\cA \rightarrow \cB$ be a normal (unital) CP-map.
A normal homomor\-phism $j\colon \cA \rightarrow \cB \otimes \cC$ 
(not necessarily unital) is a \hl{weak tensor dilation}, if there is a
normal state $\psi$ on $\cC$ such that $S = P_{\psi}\circ j$, where 
$P_{\psi}\colon \cB \otimes \cC \rightarrow \cB$ denotes the conditional
expectation determined by $P_{\psi}(b \otimes c) = b\, \psi(c)$.
\beqn{
\xymatrix{
\cA	\ar[drr]_j	\ar[rr]^S	&&	\cB						\\
					&&	\cB\otimes\cC	\ar[u]_{P_\psi=\id_\cB\otimes\psi}
}
}\eeqn
\end{definition}

\brem
Evidently, this is not a weak dilation, because the embedding of $\cB$ into 
$\cD=\cB\otimes\cC$ is unital. The possibility of a unital embedding of $\cB$ 
is due to the simple tensor product structure and should be considered as a 
benefit rather than an obstruction. We shall stay with the term weak tensor 
dilation, because it is always possible to modify $\cC$ (making it bigger) such 
that the unital embedding $\cB\rightarrow\cB\otimes\cC$ may be substituted by a 
non-unital one $\cB\rightarrow\cB\otimes p_\psi$, where $p_\psi$ is a suitable 
projection in the bigger $\cC$, such that 
$P={\eins}\otimes{p_\psi}\bullet{\eins}\otimes{p_\psi}$ does the job. 
(For instance, if the GNS-construction of $\psi$ is faithful, then we identify $\cC$ as a 
subset of $\sB(K)$ where $K$ is the GNS-Hilbert space. Then $\psi(c)=\langle k,ck\rangle$ 
for a cyclic vector $k\in K$. Enlarging $\cC$ such that the rank-one projection 
$p_\psi=|k\rangle\langle k|$ is contained, actually means setting $\cC=\sB(K)$. 
The argument may be modified suitably, if the GNS-representation is not faithful.)
In this way we can also reinterpret unital tensor dilations, extensively studied
by K\"ummerer in \cite{Kuem85}, as weak tensor dilations. Note that in 
K\"ummerer's setting, where all mappings are required unital and all states are 
required faithful, it is known that such dilations do not always exist, see 
\cite{Kuem85}, 2.1.8.
\erem

A basic result is the following existence theorem for weak
tensor dilations which confirms that weak tensor dilations
are a useful tool in the study of completely positive mappings. 

\begin{theorem}\label{mthm}
For any normal unital completely positive map $S\colon \cA \rightarrow \cB$
there exists a weak tensor dilation $j\colon \cA \rightarrow \cB \otimes \BK$,
with a vector state $\psi$ given by $k_0 \in K$.
\end{theorem}

A proof of Theorem \ref{mthm} in the case $\cA=\sB(F)$ and $\cB=\sB(G)$ using 
the Stinespring representation and the form of normal representations of 
algebras $\sB(H)$ is rather plain and well-known. 
In the general case when $\cB \subset \sB(G)$ it is not difficult
to construct in a similar way a homomorphism 
$j\colon \cA \rightarrow \BG \otimes \BK$ dilating $S$. This can even
be done on a $C^*$-algebraic level. See for example \cite{Goh04}, 1.3.3.
But the result is not a weak tensor dilation for $S\colon \cA \rightarrow 
\cB$, because the range of that $j$ need not be contained in $\cB\otimes\BK$.
 (In Section \ref{ext} we will learn that a necessary and sufficient condition 
is that a certain isometry from the space of the Stinespring representation to 
$G\otimes K$ interwines the action of the \it{commutant lifting} of $\cB'$
 and the natural action of $\cB'$ on $G\otimes K$.)

It seems that a proof of Theorem \ref{mthm}, generally, requires arguments involving von Neumann algebras. In Section \ref{proofsec} we prove Theorem \ref{mthm} by using von Neumann modules as introduced in Skeide \cite{Ske00b}. In Section \ref{ex} we illustrate the construction in a simple finite-dimensional example. After completing our proof we learned that Hensz-Chadzynska, Jajte and Paszkiewicz in \cite[Proposition 3.4]{H-CJP98} have shown a similar result with a proof based on the comparison theorem for projections in von Neumann algebras.

Our approach here underlines a far reaching duality between objects of 
categories that involve von Neumann algebras and in each case a corresponding 
category that involves the commutants of those von Neumann algebras. 
(See papers by Albeverio, Connes, Hoegh-Krohn, Muhly, Rieffel, Skeide 
and Solel \cite{Rie74a,AlHK78,Con80p,Ske03c,Ske03p,MuSo03p,MSS03p,Ske03p1}.) 
Our proof of Theorem \ref{mthm} is based on existence of a \it{complete quasi 
orthonormal system} in every von Neumann \nbd{\cB}module; see \cite{Pas73}. 
Under the duality von Neumann \nbd{\cB}modules become representations of $\cB'$ 
and existence of complete quasi orthonormal systems translates into the 
\it{amplification-induction theorem}, that is, the representation theory of 
$\cB'$; see Remark \ref{AIrem}. (Notice that in the proof  of Theorem 
\ref{mthm} we do \hl{not} apply the representation theory to the Stinespring 
representation of $\cA$. The representation to which the duality applies would 
be Arveson's \it{commutant lifting} of $\cB'$ and we do \hl{not} need the 
representation theory of the commutant lifting in the proof. The commutant 
lifting has its appearance not before Section \ref{ext}.)

A further instance of the duality is the correspondence (see Albeverio and Hoegh-Krohn \cite{AlHK78}) between unital CP-maps $S\colon\cA\rightarrow\cB$ and $S'\colon\cB'\rightarrow\cA'$ where the maps are covariant with respect to certain vector states. In Theorem \ref{dualthm} we reprove a duality from \cite{Goh04,Goh03a} between covariant extensions $Z\colon \sB(F) \rightarrow \sB(G)$ of $S$ and weak tensor dilations of $S'$. This makes Theorem \ref{mthm} actually two theorems, one on existence of weak tensor dilations (Theorem \ref{mthm}) and another one on existence of extensions of CP-maps (Theorem \ref{extexprop}). As explained in \cite{Goh04}, there are applications in the theory of noncommutative Markov processes and in particular in K\"ummerer-Maassen scattering theory \cite{KueMa00}. The statement of this duality has been the original motivation for the notion of weak tensor dilation because here the relaxed versions of dilation mentioned above are not sufficient. By applying the methods of \cite{Ske03p1} for a reconstruction of the
correspondence between dilations and extensions we put not only that duality 
but a couple of others into a new unified perspective. Combining the correspondence with the existence of weak tensor dilations given by Theorem \ref{mthm}, we infer an existence result for extensions of completely positive maps with states (Theorem \ref{extexprop}), which is a refinement of a well known extension result of Arveson \cite{Arv69}.

We illustrate the connections in the following diagram.
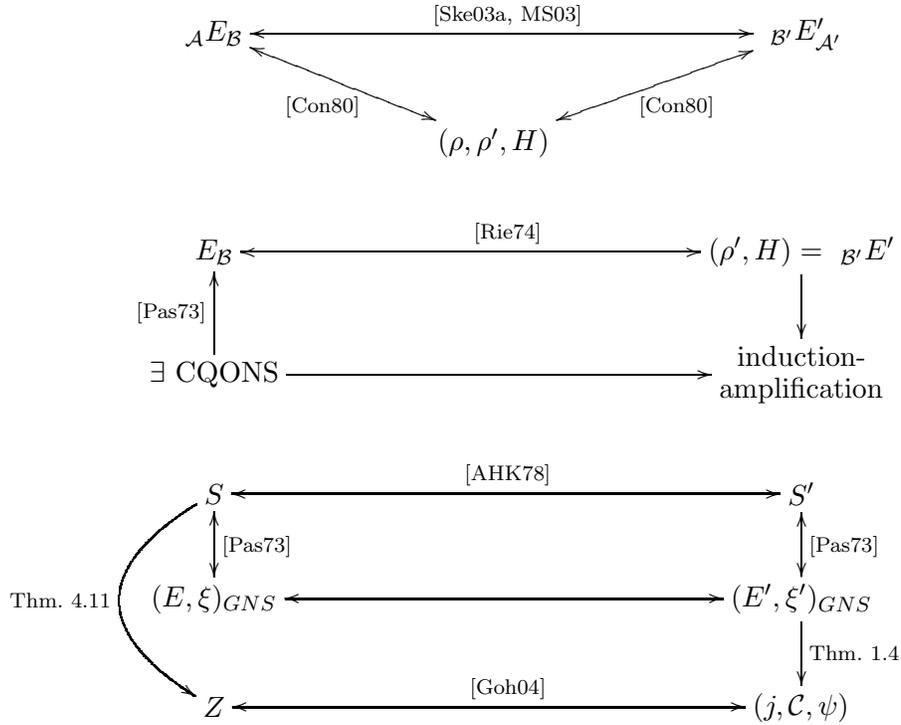
\beqn{
\xymatrix{
_\cA E_\cB	\ar@{<->}[rrrr]^{\text{\cite{Ske03c,MuSo03p}}}	&&&&{~}_{\cB'}E'_{\cA'}
\\
&&	(\rho,\rho',H)		\ar@{<->}[ull]^{\text{\cite{Con80p}}~~}	\ar@{<->}[urr]_{~~\text{\cite{Con80p}}}	&&
\\
E_\cB		\ar@{<->}[rrrr]^{\text{\cite{Rie74a}}}	&&&&(\rho',H)={~}_{\cB'}E'	\ar[d]
\\
\txt{$\exists$ CQONS}	\ar[u]^{\text{\cite{Pas73}}}	\ar[rrrr]	&&&&\txt{induction-\\amplification}
\\
S	\ar@/_3pc/[dd]_{\text{Thm.\ \ref{extexprop}}}			\ar@{<->}[d]^{\text{\cite{Pas73}}}		\ar@{<->}[rrrr]^{\text{\cite{AlHK78}}}	&&&&	S'	\ar@{<->}[d]^{\text{\cite{Pas73}}}		
\\
(E,\xi)_{GNS}	\ar@{<->}[rrrr]				&&&&	(E',\xi')_{GNS}	\ar[d]^{\text{Thm.\ \ref{mthm}}}
\\
Z	\ar@{<->}[rrrr]^{\text{\cite{Goh04}}}	&&&&	(j,\cC,\psi)
}
}\eeqn
Here $\rho$ and $\rho'$ are a pair of representations of $\cA$ and of $\cB'$, respectively, on the same Hilbert space whose ranges mutually commute. $(E,\xi)_{GNS}$ and $(E',\xi')_{GNS}$ denote the GNS-modules and cyclic vector of $S$ and $S'$, respectively. In the case of these GNS-constructions $\rho$ is the Stinespring representation of $\cA$ for $S$ and the commutant lifting of $\cA$ for $S'$, while $\rho'$ is the commutant lifting of $\cB'$ for $S$ and the Stinespring representation of $\cB'$ for $S'$. (We see: Looking only at the Stinespring representation of a CP-map misses half the information. Only taking also into account Arveson's commutant lifting gives full information. However, looking at the pair is the same as looking immediately at the GNS-module.)

\section{Proof of Theorem \ref{mthm}}\label{proofsec}

Von Neumann modules can be characterized as self-dual Hilbert modules over von Neumann 
algebras. In this algebraic and intrinsic form the characterization is suitable also for 
modules over $W^*$--algebras. For these our essential tool here, existence of 
\it{complete quasi orthonormal systems}, was proved already by Paschke \cite{Pas73}. 
We put, however, emphasis on considering a von Neumann algebra $\cB$ acting 
(always non-degenerately) as a concrete strongly closed subalgebra of operators 
\it{on} a Hilbert space $G$ and, following Skeide \cite{Ske00b}, we consider von Neumann 
$\cB$--modules $E$ as concrete strongly closed submodules of operators \it{between} Hilbert 
spaces $G$ and $H$. Once the identifying representation of $\cB$ on $G$ is fixed, the 
construction of $H$ and of the embedding $E\rightarrow\sB(G,H)$ is canonical. 

We start by recalling the definition and basic properties of von Neumann modules as 
introduced in \cite{Ske00b}. For a more detailed account we refer the reader to Skeide 
\cite{Ske01,Ske03p1}.

Let $\cB\subset\sB(G)$ be a von Neumann algebra. Furthermore, let $E$ be a Hilbert $\cB$--module. Then the (interior) tensor product (over $\cB$) $H=E\odot G$ of the Hilbert $\cB$--module $E$ and the Hilbert $\cB$--$\C$--module $G$ is a Hilbert $\C$--module, i.e.\ a Hilbert space, with inner product $\AB{x\odot g,x'\odot g'}=\AB{g,\AB{x,x'}g'}$. Every element $x\in E$ gives rise to a mapping $L_x\in\sB(G,H)$ defined by $L_xg=x\odot g$ with adjoint $L_x^*(y\odot g)=\AB{x,y}g$. One verifies that $\AB{x,y}=L_x^*L_y$ and that $L_{xb}=L_xb$. In other words, we may and, from now on, we will identify a Hilbert module over a von Neumann algebra $\cB\subset\sB(G)$ as a submodule of $\sB(G,H)$ with the natural operations.

An (adjointable and, therefore, bounded) operator $a\in\sB^a(E)$ gives rise to an operator $x\odot g \mapsto ax\odot g$ in $\sB(H)$. We, therefore, may and will identify the $C^*$--algebra $\sB^a(E)$ as a subalgebra of $\sB(H)$.

\bdefi
A Hilbert module $E$ over a von Neumann algebra $\cB\subset\sB(G)$ is a \hl{von Neumann $\cB$--module}, if $E$ is strongly closed in $\sB(G,H)$.
\edefi

It is immediate that in this case $\sB^a(E)$ is a von Neumann subalgebra of $\sB(H)$.

The basic tool in establishing the criterion from \cite{Ske00b,Ske01} that $E$ 
is a von Neumann module, if and only if it is \hl{self-dual} (i.e.\ every 
bounded right linear mapping $E\rightarrow\cB$ has the form $\AB{x,\bullet}$ 
for a unique $x\in E$) was the construction of a quasi orthonormal system. 
(See, however, Skeide \cite{Ske03p} for a new quick proof without using quasi orthonormal systems based on methods from \cite{Ske03p1}.)

\bdefi
A \hl{quasi orthonormal system} is a family $\bfam{e_i,p_i}_{i\in I}$ of pairs consisting of an element $e_i\in E$ and a non-zero projection $p_i\in\cB$ such that
\beqn{
\langle e_i,e_{i'}\rangle
=
p_i\delta_{ii'}.
}\eeqn
We say the family is \hl{orthonormal}, if $p_i=\U$ for all $i\in I$.
\edefi

In particular, the $e_i$ are partial isometries and the net of projections
\beqn{
\Bfam{\sum_{i\in I'}e_i e_i^*}_{I'\subset I,\#I'<\infty}
}\eeqn
is increasing in $\sB^a(E)$ over the finite subsets $I'$ of $I$ and converges, therefore, 
to an element $p_I\in\sB^a(E)$. Obviously, $p_I$ is the projection onto the von Neumann 
$\cB$--submodule of $E$ generated by $\bfam{e_i}_{i\in I}$.

\bdefi
A quasi orthonormal system $\bfam{e_i,p_i}_{i\in I}$ is called \hl{complete}, if $p_I=\eins$.
\edefi

Of course, $\bfam{e_i,p_i}_{i\in I}$ is a complete quasi orthonormal system for $p_I E$. Existence of a complete quasi orthonormal system follows now as a standard application of Zorn's lemma. In particular, a complete quasi orthonormal system  can be chosen such that it contains a given quasi orthonormal system. 

\bprop\label{embed}
Every von Neumann $\cB$--module is isomorphic to the (von Neumann module) direct sum
\beqn{
E
~=~\bigoplus_{i\in I}p_i\cB
}\eeqn
of right ideals $p_i\cB$ for some complete quasi orthonormal system $\bfam{e_i,p_i}_{i\in I}$. It is, therefore, isomorphic to a complemented submodule of
\beqn{
\bigoplus_{i\in I}\cB
~=~
\cB\otimes K}\eeqn
where $K$ is a Hilbert space with orthonormal basis $\bfam{k_i}_{i\in I}$ and $\cB\otimes K$ 
is the strong closure of the (right) $\cB$--submodule of $\sB(G,G\otimes K)$ generated by 
the mappings $g\mapsto g\otimes k$ ($k\in K)$.
\eprop

Explicitly, an element $e_i b \in E$ corresponds to $p_i b \otimes k_i \in \cB\otimes K$ and is represented by the mapping $g\mapsto p_i b g\otimes k_i$ in $\sB(G,G\otimes K)$. Clearly, $\bfam{e_i,p_i}_{i\in I}$ is still a quasi orthonormal system (in general, not complete, of course) for $\cB\otimes K$ and $p_I$ is the projection onto $E$.

\bdefi
Let $\cA$ be another von Neumann algebra. A von Neumann $\cA$--$\cB$--module is an $\cA$--$\cB$--module and a von Neumann $\cB$--module such that the left action defines a normal (unital, because the action of $\cA$ is unital) representation of $\cA$ on $H$.
\edefi

\brem\label{Stinerem}
Let $S\colon\cA\rightarrow\cB$ be a normal CP-map. Then the strong closure of the GNS-module $E$ of $S$ (that is the unique Hilbert $\cA$--$\cB$--module generated by a single element $\xi$ which fulfills $\AB{\xi,a\xi}=S(a)$; see Paschke \cite{Pas73}) is a von Neumann $\cA$--$\cB$--module; see Bhat and Skeide \cite{BhSk00} and \cite{Ske01}.

Observe that $\cA\rightarrow\sB(H)$ is the Stinespring representation and that $\xi^*a\xi=S(a)$.
\erem

\lf\noindent
\sc{Proof of Theorem \ref{mthm}.~} Let $S\colon\cA\rightarrow\cB$ be a unital CP-map and 
let $E$ be its GNS-module with cyclic vector $\xi$. Since $S$ is unital, $\xi$ is a unit 
vector, i.e.\ $\AB{\xi,\xi}=\eins$. Let $\bfam{e_i,p_i}_{i\in I}$ be a complete quasi 
orthonormal system  which contains $\xi=e_0$ and identify $E$ as a submodule of 
$\cB\otimes K$ as described in Proposition \ref{embed}. Then $\xi \simeq \eins \otimes k_0$.

Now $j\colon a\mapsto a p_I$ defines a (usually non-unital) normal homomorphism 
$\cA\rightarrow\sB^a(\cB\otimes K)$. There is a well-known identification 
$\sB^a(\cB\otimes K) = \cB\otimes\sB(K)$, see \cite{Ske01}. Explicitly, 
$|b_1\otimes\eta_1\rangle\langle b_2\otimes\eta_2| \in \sB^a(\cB\otimes K)$ is identified 
with $b_1 b^\star_2 \otimes |\eta_1\rangle\langle\eta_2| \in \cB\otimes\sB(K)$. Then we can 
compute
\beqn{
(\id_\cB\otimes\AB{k_0,\bullet k_0})\circ j(a)
~=~
\AB{\eins \otimes k_0, j(a) \eins \otimes k_0} = \AB{\xi, a \xi} = S(a).
}\eeqn
In other words, $(j,\AB{k_0,\bullet k_0})$ is a weak tensor dilation of $S$ 
(to $\cB\otimes\sB(K)$).\qedsymbol

\brem
The result can be saved partly, if $S$ is not unital. Then $\xi$ is not a unit vector, but it posesses a \hl{polar decomposition} $\xi_0\abs{\xi}$ with $\abs{\xi}=\sqrt{\AB{\xi,\xi}}$. In this case, $\AB{\xi_0,\xi_0}=p_0$ is a projection. We do the same construction and may express $S$ as $\bfam{(\abs{\xi}\bullet\abs{\xi})\otimes\AB{k_0,\bullet k_0}}\circ j$.
\erem

\section{An example}\label{ex}

The following example illustrates the construction
of the preceding section. We consider the stochastic matrix
$S=\frac{1}{2}\bfam{\substack{1~1\\1~1}}$,
i.e. $\cA = \cB = \C^2$ and 
$S\bfam{\substack{a_1\\a_2}}=\frac{1}{2}\bfam{\substack{a_1 + a_2 \\a_1 + a_2 }}$. 
Then the GNS-module is $E = \cA \otimes \cB = \C^2 \otimes \C^2$. In fact, for $x_1,x_2,y_1,y_2 \in \cA$ and $p_1 = \bfam{\substack{1\,\\0}},p_2 = \bfam{\substack{0\\1\,}}$ the inner product for $x = x_1 \otimes p_1 + x_2 \otimes p_2$ and $y = y_1 \otimes p_1 + y_2 \otimes p_2$ takes the form
\beqn{
\AB{x,y}
~=~
\sum^2_{i,j=1} p^\star_i S(x^\star_i y_j) p_j= p_1\, S(x^\star_1\, y_1) + p_2\, S(x^\star_2\, y_2).
}\eeqn
In particular, there are no non-trivial vectors of zero length. With this formula we can also check 
that
\[
\bCB{e_0 := \eins \otimes \eins,\quad e_1 := \bfam{\substack{1\,\\-1}}\otimes p_1,\quad e_2 := \bfam{\substack{1\,\\-1}}\otimes p_2}
\]
is a complete quasi orthonormal system for $E$. In fact, its $\cB$-linear span is $E$ and we have
\[
\langle e_0,e_0 \rangle = \eins,~~\langle e_1,e_1 \rangle = p_1,~~\langle e_2,e_2 \rangle = p_2 \quad\mbox{and}\quad\langle e_i,e_j \rangle = 0\; \mbox{for}\; i \not= j.
\]
Now let $K$ be a $3$-dimensional Hilbert space with ONB $\{k_0, k_1, k_2\}$ as in 
Proposition \ref{embed}, so that $E$ can be identified with the module $p_I(\cB \otimes K)$ for
\[
p_I
~=~
\Bfam{\substack{\eins~~~~~~\\~~~p_1~~~\\~~~~~~p_2}}
~\in~
M_3(\C^2)
~\simeq~
\cB \otimes \BK.
\]
According to Theorem \ref{mthm} and its proof, we can now construct a weak tensor dilation $j\colon \cA \rightarrow \cB \otimes \BK$ by $j(a) := a\, p_I$, where we use the left action of $\cA$ on $E\subset\cB\otimes K$. Let us compute explicitly the matrix of $j(a)$ for $a =\bfam{\substack{a_1\\a_2}}\in \cA$ with respect to the $\cB$--module basis $\{\eins \otimes k_0, \eins \otimes k_1, \eins \otimes k_2\}$ of $\cB \otimes \cK$. With
\[
a\, e_0
~=~
\bfam{\substack{a_1\\a_2}}\otimes \eins
~=~
{\textstyle\frac{1}{2}}(a_1+a_2) e_0 + {\textstyle\frac{1}{2}}(a_1-a_2) (e_1+e_2)
\]
and with $p_I \eins \otimes k_i = p_i \otimes k_i \simeq e_i$ we find that
\[
j(a) \eins \otimes k_0
~=~
{\textstyle\frac{1}{2}}(a_1+a_2)\, \eins \otimes k_0 + {\textstyle\frac{1}{2}}(a_1-a_2)\, (p_1 \otimes k_1 + p_2 \otimes k_2).
\]
Thus, for example,
\[
j(a)_{10}
~=~
\AB{\eins \otimes k_1, j(a) \eins \otimes k_0}
~=~
{\textstyle\frac{1}{2}}(a_1-a_2) p_1.
\]
By similar computations we get
\begin{equation*}
\begin{split}
& j(a)
~\simeq~
\frac{1}{2}\left( \begin{array}{ccc}        a_1+a_2           & (a_1-a_2) p_1 & (a_1-a_2) p_2 \\         (a_1-a_2) p_1 & (a_1+a_2) p_1 &    0              \\         (a_1-a_2) p_2 &    0              & (a_1+a_2) p_2 \\        \end{array} \right)
\\
\simeq~&p_1 \otimes \frac{1}{2}\left( \begin{array}{ccc}        a_1+a_2 & a_1-a_2 & 0 \\         a_1-a_2 & a_1+a_2 & 0 \\           0     &   0     & 0 \\        \end{array} \right) +p_2 \otimes \frac{1}{2}\left( \begin{array}{ccc}        a_1+a_2 & 0 & a_1-a_2 \\           0     & 0 &   0     \\         a_1-a_2 & 0 & a_1+a_2 \\        \end{array} \right)
\end{split}
\end{equation*}
in $\cB \otimes \BK$. Using different methods, in \cite{Goh04,Goh03a} 
a complete classification of weak tensor dilations for this particular map $S$ has been 
obtained. Comparing the matrix for $j(a)$ with this classification, it is easy to check 
that the weak tensor dilation $j$ constructed above coincides with $j_0$ in \cite[1.7.1]
{Goh04}.

\section{Weak tensor dilations and extensions of CP-maps revisited}\label{ext}

In this section we want to reconstruct the duality theory for dilations and extensions of \cite{Goh04,Goh03a} in terms of von Neumann modules. In the end, we will deal with a pair of CP-maps $S\colon\cA\rightarrow\cB$ and $S'\colon\cB'\rightarrow\cA'$ which, under certain conditions, are related by a duality. Extensions to mappings $\sB(F)\rightarrow\sB(G)$ of $S$ correspond, then, to weak tensor dilations of $S'$. For the sake of clarity, however, we introduce at each step only as much structure as is necessary to make the argument work.

Let $\cB\subset\sB(G)$ be a von Neumann algebra and let $E$ be a von Neumann \nbd{\cB}module. As in Section \ref{proofsec}, we identify $E$ as a subset of $\sB(G,H)$ with $H=E\odot G$ and $xg=x\odot g$. A new ingredient comes from the observation that there is an associated normal unital representation $\rho'$ of the commutant $\cB'$ of $\cB\subset\sB(G)$ on $H=E\odot G$ where $\rho'(b')$ lets act an element $b'\in\cB'$ in the natural way on $x\odot g$ as $\id_E\odot b'$. (Since $b'$ is a bilinear operator on the $\cB$--$\C$--module $G$, the operator $\rho'(b')$ is well-defined, and it is routine to check that $\rho'$ is normal.)

\brem
In the case when $E$ is the (strong closure of the) GNS-module of a CP-map (see Remark \ref{Stinerem}), thus, $H$ is the representation space of the Stinespring construction, this representation of $\cB'$ is known as \hl{commutant lifting} from Arveson \cite{Arv69}.
\erem

We recover $E$ as the intertwiner space
\beqn{
C_{\cB'}(\sB(G,H))
~=~
\bCB{x\in\sB(G,H)\colon \rho'(b')x=xb'~(b'\in\cB')}.
}\eeqn
To see this, we observe that $C_{\cB'}(\sB(G,H))$ is a von Neumann $\cB$--module which, clearly, contains $E$. On the other hand, the complement of $E$ is $\CB{0}$, since $\cls EG=H$. In analogy with Hilbert spaces, a submodule of a von Neumann module with zero-complement is strongly total (see \cite{Ske00b,Ske01}), and since $E$ is strongly closed it must be all of $C_{\cB'}(\sB(G,H))$. (In other words, for a given representation $\rho'$ of $\cB'$ on $H$ there is at most one von Neumann \nbd{\cB}submodule $E$ of $\sB(G,H)$ such that $H=E\odot G=\cls EG$ and $\rho'$ is the commutant lifting.)

The converse of the observation that every normal unital representation $\rho'$ of $\cB'$ on a Hilbert space $H$ gives rise to a von Neumann $\cB$--module $E=C_{\cB'}(\sB(G,H))$ can be traced back already to Rieffel \cite{Rie74a}. (The crucial point for the direction here is a lemma from Muhly and Solel \cite{MuSo02} which asserts that $\cls C_{\cB'}(\sB(G,H))G=H$.) In other words, von Neumann $\cB$--modules and normal unital representations of $\cB'$ are two ways to speak about the same object; see \cite{Ske03p1}.

\brem\label{com}
It is not difficult to see that $\sB^a(E)$ is the commutant of $\rho'(\cB')$ in $\sB(H)$. In fact, the correspondence 
$E\subset\sB(G,H)\leftrightarrow(\rho',H)$ extends to a functor which sends the morphisms $a\in\sB^a(E_1,E_2)$ to the morphisms $a\odot\id_G\in\sB^{bil}(H_1,H_2)$ and back.
\erem

\brem\label{AIrem}
The so-called amplification-induction theorem, which summarizes the representation theory of 
von Neumann algebras, (see, e.g., \cite[IV.5.5]{Tak79}) 
is now a corollary of Proposition \ref{embed}. Concretely, 
given a normal unital representation 
$\rho'$ of $\cB'$ on $H$, there is a Hilbert space $K$ and a projection 
$p_I\in(\cB'\otimes\id_K)'=\cB\otimes\sB(K)$ such that 
$H=p_I(G\otimes K)$ and $\rho'(b')=p_I(b'\otimes\id_K)$. Additionally, we can say that 
$p_I=\bigoplus_{i\in I}p_i\otimes|k_i\rangle\langle k_i|$ is diagonal in the basis of $K$.
\erem

Now suppose $E$ is a von Neumann \nbd{\cA}\nbd{\cB}module where $\cA\subset\sB(F)$ is another von Neumann algebra. On $H=E\odot G$ we have the Stinespring representation $\rho$ of $\cA$ and the commutant lifting $\rho'$ of $\cB'$. By the preceding discussion we can recover the right module $E$ as intertwiner subspace of $\sB(G,H)$ for $\cB'$ via $\rho'$ and its left action via $\rho$. In other words, von Neumann $\cA$--$\cB$--modules and pairs of normal unital representations $\rho$ of $\cA$ and $\rho'$ of $\cB'$ on the same Hilbert space $H$ with commuting range ($\SB{\rho(\cA),\rho'(\cB')}=\CB{0}$) are two sides of  the same coin.

By an observation from Skeide \cite{Ske03c} (for the case $\cA=\cB$) and Muhly and Solel \cite{MuSo03p} (for the general case) nobody prevents us from exchanging the roles of $\cA$ and $\cB'$. In other words, we can construct the intertwiner space
\beqn{
E'
~=~
C_\cA(\sB(F,H))
~=~
\bCB{x'\in\sB(F,H))\colon\rho(a)x'=x'a~(a\in\cA)}
}\eeqn
which is a von Neumann $\cA'$--module and which inherits a left action of $\cB'$ via $\rho'$. In other words, there is a duality between von Neumann $\cA$--$\cB$--modules and von Neumann $\cB'$--$\cA'$--modules. We call $E'$ the \hl{commutant} of $E$. (The notion of commutant of $\cB$ is contained when $E$ is the von Neumann $\cB$--$\cB$--module $\cB$.) Clearly, $E''=E$ and there is also the analogue of von Neumann's double commutant theorem.

\brem\label{MSSrem}
In \cite{Ske03p1} this duality plays a crucial role in developing the complete theory of normal representations of $\sB^a(E)$ for a von Neumann $\cB$--module by adjointable operators on a von Neumann $\cC$--module generalizing Arveson's technique of intertwiner spaces for $\sB(G)$. Meanwhile, in Muhly, Skeide and Solel \cite{MSS03p} there exists a simpler approach that works already in the \nbd{C^*}case (of course, without nice relations to commutants).
\erem

So far, we have a duality, the commutant, between von Neumann $\cA$--$\cB$--modules $E\subset\sB(G,H)$ and von Neumann $\cB'$--$\cA'$--modules $E'\subset\sB(F,H)$. Now suppose that $S$ is a normal unital CP-map from a von Neumann algebra $\cA\subset\sB(F)$ to a von Neumann algebra $\cB\subset\sB(G)$. Furthermore, let $E$ be a von Neumann \nbd{\cA}\nbd{\cB}module that contains a vector $\xi$ such that $S(a)=\AB{\xi,a\xi}$. (That is, the von Neumann \nbd{\cA}\nbd{\cB}submodule of $E$ generated by $\xi$ is the (strong closure of the) GNS-module of $S$.) We ask, whether we can extend the duality between the modules $E$ and $E'$ to a duality between (unital normal) CP-maps $S\colon\cA\rightarrow\cB$ and $S'\colon\cB'\rightarrow\cA'$. In other words, given $E$ with the unit vector $\xi$, can we construct in a ``reasonable'' way a unit vector $\xi'\in E'$ (defining a unital normal CP-map $S'=\AB{\xi',\bullet\xi'}$)? Can we apply the construction to $\xi'\in E'$ such that it gives back $\xi$ completing the duality?

Such a construction is possible but needs for additional structure in either direction. Let us start with the assumptions for the construction of $\xi'\in E'$ from $\xi\in E$. So, as before, let $E$ be a von Neumann \nbd{\cA}\nbd{\cB}module with commutant $E'$ and with a unit vector $\xi\in E$ so that $S=\AB{\xi,\bullet\xi}$ defines a normal unital CP-map $\cA\rightarrow\cB$. Furthermore, suppose that $f$ and $g$ are unit vectors in $F$ and $G$ such that the vector states $\vp_f=\AB{f,\bullet f}$ and $\vp_g=\AB{g,\bullet g}$ are \hl{covariant} for $S$, i.e.\ $\vp_f=\vp_g\circ S$. Finally, suppose that $f$ is cyclic for $\cA$. Then, with the help of the vector $\xi \in E$, we can define an isometry $\xi'\colon F\rightarrow H$by setting
\beq{\label{xi'def}
\xi'\colon af
~\longmapsto~a\xi\odot g.
}\eeq
By definition $\xi'$ intertwines the actions of $\cA$ and, therefore, is a unit vector in $E'$. Using $\rho'(b')\xi=\xi b'$, we obtain the following result.

\bprop
The mapping $S'\colon\cB'\rightarrow\cA'$ defined by $S'(b')=\AB{\xi',b'\xi'}$ is the unique (unital normal) CP-map fulfilling $\vp_f\circ S'=\vp_g$ and $\vp_g(b'S(a))=\vp_f(S'(b')a)$. It is called the \hl{dual map} of $S$.

Moreover, exchanging the roles of $\cA$ and $\cB'$, if also $g$ is cyclic for $\cB'$, then $S''=S$ so that we obtain a true duality $S\leftrightarrow S'$.
\eprop

\brem
For faithful states and from a different point of view this duality was discussed by Albeverio and Hoegh-Krohn in \cite{AlHK78}.
\erem

\brem\label{GNSdual}
If, under the cyclicity conditions for $f$ and $g$, the vector $\xi$ is cyclic for $E$, i.e.\ if $E$ is the GNS-module of $S$, then $\xi'$ is cyclic for $E'$, i.e.\ $E'$ is the GNS-module of $S'$. Indeed, $\xi$ is cyclic for $E$, if and only if $\cls\cA\xi G=H$ (because then the orthogonal complement is $\zero$). Likewise, $\xi'$ is cyclic for $E'$, if and only if $\cls\cB'\xi'F=H$. By \eqref{xi'def} we find
\beqn{
b'\xi'af
~=~
a\xi b'g,
}\eeqn
so that cyclicity of $f$ for $\cA$, cyclicity of $\xi$ for $E$ and cyclicity of $g$ for $\cB'$ imply cyclicity of $\xi'$ for $E'$.
\erem

\brem
What happens, if the conditions are weakened? If, for instance, $S$ is not unital, i.e.\ if $\xi$ is not a unit vector but $S$ still satisfies the covariance condition, then \eqref{xi'def} still defines an isometry. That is for the backwards direction ($g$ assumed cyclic for $\cB'$) we will never get back $\xi$. On the other hand, if $f$ is not cyclic for $\cA$ we can still use \eqref{xi'def} to define a (unique) partial isometry with cokernel $\ol{\cA f}$ resulting in a possibly non-unital CP-map $S'$. Also we may, in general, not hope to recover $\xi$ without $g$ being cyclic for $\cB'$. So, partly, the constructions work also under more general assumptions, but they do not yield as nice structures.
\erem

Now we can construct the duality between weak tensor dilations of $S'$ and 
extensions of $S$ for certain covariant vector states $\vp_f$ and $\vp_g$. 
Every direction involves the construction of that mapping ($S$ or $S'$) that 
is not given in the beginning. Correspondingly, every direction needs only 
one of the cyclicity conditions for the duality between $S$ and $S'$.

Let us start with a weak tensor dilation $(j\colon\cB'\rightarrow\cA'\otimes\sB(L),\vp_\ell=\AB{\ell,\bullet\ell})$ of a given CP-map $S'\colon\cB'\rightarrow\cA'$ where $L$ is  Hilbert space and $\ell$ is a unit vector in $L$. Then $\cA'\otimes L\subset\sB(F,F\otimes L)$ as a von Neumann \nbd{\cA'}module is exactly the intertwiner subspace of $\sB(F,F\otimes L)$ for the representation $a\mapsto a\otimes\id_L$ of $\cA$. Therefore, by Remark \ref{com}, $\cA'\otimes\sB(L)=(\cA\otimes\id_L)'$ is exactly the algebra $\sB^a(\cA'\otimes L)$ of adjointable operators on 
$\cA'\otimes L$. In particular, $j(\eins_{\cB'})$ is a projection in $\sB^a(\cA'\otimes L)$.
Furthermore,
\beqn{
(\id_F\otimes\ell)^*j(b')(\id_F\otimes\ell)
~=~
S'(b')
}\eeqn
so that $E'=j(\eins_{\cB'})(\cA'\otimes L)$ is a von Neumann \nbd{\cB'}\nbd{\cA'}submodule of $\cA'\otimes L$ with left 
multiplication via $j$ which contains a vector $\xi'=(\id_F\otimes\ell)\in E'\subset\cA'\otimes L$ such that $\xi'^*b'\xi'=S'(b')$. In other words, we have identified $H=E'\odot F$ with the subspace $j(\eins_{\cB'})(F \otimes L)$ of $F \otimes L$. Because $a\otimes\id_L$ commutes with $j(\eins_{\cB'})$ it restricts to the operator $(a\otimes\id_L)j(\eins_{\cB'})=\rho(a)$ on $H$. We are now in the situation as described after Remark \ref{MSSrem}. In particular, we have the von Neumann \nbd{\cA}\nbd{\cB}module
\beqn{
E
~=~
E''
~=~
C_{\cB'}(\sB(G,H))
~\subset~
\sB(G,F\otimes L)
}\eeqn
(with left multiplication via $\rho$) and, supposing that $g$ is cyclic for $\cB'$, we get 
an isometry $\xi\colon G\rightarrow H\subset F\otimes L$ in $E$ by setting
\beqn{
\xi(b'g)
~=~
b' \xi'\odot f
~=~
j(b')(f\otimes\ell)
}\eeqn
so that
\beqn{
S(a)
~=~
\xi^*\rho(a)\xi
~=~
\xi^*(a\otimes\id_L)\xi
}\eeqn
is the dual map $S=S''$ of $S'$. (In \cite{Goh04,Goh03a} $\xi$ is called the isometry \hl{associated} with the dilation $j$ and with $g$.) Summarizing, the weak tensor dilation of $S'$ provides us with a larger space 
$F \otimes L$ and with corresponding embeddings of our objects and cyclicity of $g$ for $\cB'$ provides us with the dual CP-map $S$ of $S'$. The following observation is now easily checked: The mapping $Z\colon x\mapsto\xi^*(x\otimes\id_L)\xi$ $(x\in\sB(F))$ is an extension of $S$ to a unital 
completely positive mapping $\sB(F)\rightarrow\sB(G)$ which fulfills $\vp_f=\vp_g\circ Z$.

Conversely, let $Z\colon\sB(F)\rightarrow\sB(G)$ be a unital normal completely positive extension of a given CP-map $S\colon\cA\rightarrow\cB$ fulfilling $\vp_f=\vp_g\circ Z$. Like every unital completely positive mapping $\sB(F)\rightarrow\sB(G)$, also $Z$ can be written in the form $x\mapsto\xi^*(x\otimes\id_L)\xi$ for a suitable (isometric) mapping $\xi\in\sB(G,F\otimes L)$. (This is so, because all von Neumann $\sB(F)$--$\sB(G)$--modules have the particularly simple form $\sB(G,F\otimes L)$ for a suitable Hilbert space $L$ with the canonical actions of $\sB(F)$ and $\sB(G)$; see \cite{BhSk00,Ske01}.) Let $E$ denote the strongly closed $\cA$--$\cB$--submodule of $\sB(G,F\otimes L)$ generated by $\xi$. Clearly, $E$ is the GNS-module of $S$ and $\xi$ is the cyclic vector. More precisely, $H=\cls(\cA\otimes\id_L)\xi\cB G=E\odot G$ and $\xi=p_H\xi$. Therefore, on $H$ there is the commutant lifting $\rho'$ of $\cB'$ which we can use to define a (usually non-unital) homomorphism $j(b')=\rho'(b')p_H$. It is clear that $j(b')$ and $a\otimes\id_L$ commute, because $(a\otimes\id_L)p_H\in\sB^a(E)$; cf.\ Remark \ref{com}. Furthermore,
\beqn{
\AB{f,xf}
~=~
\AB{g,\xi^*(x\otimes\id_L)\xi g}
~=~
\AB{\xi g,(x\otimes\id_L)\xi g},
}\eeqn
and, in particular,
\beqn{
1
~=~
\bAB{f,\bfam{|f\rangle\langle f|}f}
~=~
\bAB{\xi g,\bfam{|f\rangle\langle f|\otimes\id_L}\xi g}.
}\eeqn
Therefore, $\xi g$ is in the range $f\otimes L$ of the projection $|f\rangle\langle f|\otimes\id_L$ and $\ell=(\langle f|\otimes\id_L)\xi g$ is the unique unit vector in $L$ such that $\xi g=f\otimes\ell$. Define $\vp_\ell=\AB{\ell,\bullet\ell}$. Supposing now that $f$ is cyclic for $\cA$,  we define $S'$ as the dual CP-map of $S$. We find
\bmun{
\AB{a_1f,(\id_F\otimes\vp_\ell)\circ j(b')a_2f}
~=~
\AB{f\otimes\ell,(a_1^*\otimes\id_L)j(b')(a_2\otimes\id_L)(f\otimes\ell)}
\\
~=~
\AB{\xi g,(a_1^*a_2\otimes\id_L)j(b')\xi g}
~=~
\AB{\xi g,(a_1^*a_2\otimes\id_L)p_H\xi b'g}
~=~
\AB{g,S(a_1^*a_2)b'g}
\\
~=~
\AB{f,a_1^*a_2S'(b')f}
~=~
\AB{a_1f,S'(b)a_2f}.
}\emun
Since $f$ is cyclic for $\cA$, it follows that $(\id_F\otimes\vp_\ell)\circ j(b')=S'(b')$. In other words, $(j,\vp_\ell)$ is a 
weak tensor dilation of $S'$.

\brem
While in the direction $Z\to j$ ($f$ cyclic) the module $E$ is the GNS-module of $S$ (hence, if $g$ is also cyclic, then by Remark \ref{GNSdual} $E'$ is the GNS-module of $S'$, too), in the opposite direction $j\to Z$ ($g$ cyclic) the situation is more general in that here $E'$ need only to contain the GNS-module of $S'$. If $E'$ is the GNS-module of $S'$ (and if $g$ and $f$ are cyclic), then it is easy to check that the construction $Z\to j$ applied to the result of $j\to Z$ gives back the same weak dilation we started with. The difference between $E'\supset\mrm{GNS}_{S'}$ and $E'=\mrm{GNS}_{S'}$ is the same as the difference between an arbitrary weak tensor dilation and its minimal version. Under a suitable notion of equivalence for weak tensor dilations (where every weak tensor dilation is equivalent to its unique minimal version) the correspondence between extensions and dilations is, indeed, one to one. The relation between this subtle  equivalence relation among weak tensor dilations (finer than just unitary equivalence, because all minimal weak dilations are unitarily equivalent) and a certain ``subconvex'' structure among the expansion coefficients with respect to some QONB (\it{Kraus decomposition} of $S$) will be investigated elsewhere.

These correspondences between completely positive extensions of $S$ and weak tensor dilations of $S'$ (under the conditions that $f$ is cyclic for $\cA$ for the one direction and that $g$ is cyclic for $\cB'$ for the other direction) are the main results of the first chapter of \cite{Goh04}. We recover them here by more module orientated methods and put them into perspective with the recent 
developments in \cite{Ske03p1} as a part of a more far reaching duality between ``von Neumann objects'' and their \it{commutants}.
\erem

Let us summarize the duality in form of a theorem (due to \cite{Goh04}) and then combine it with Theorem \ref{mthm}.

\begin{theorem}\label{dualthm}
Let $\cA\subset\sB(F)$ and $\cB\subset\sB(G)$ be von Neumann algebras with 
vector states $\vp_f$ $(f\in F)$ and $\vp_g$ $(g\in G)$ where $f$ is cyclic 
for $\cA$ and $g$ is cyclic for $\cB'$. Then the duality $S \leftrightarrow S'$ 
between CP-maps for which the states are covariant, extends to a
(still one-to-one) duality between normal unital completely positive covariant 
extensions $Z\colon\sB(F)\rightarrow\sB(G)$ of $S$ and minimal weak tensor 
dilations of $S'$.
\end{theorem}

\begin{theorem}\label{extexprop}
Suppose that $S$ is a normal unital CP-map from a von Neumann algebra 
$\cA\subset\sB(F)$ to a von Neumann algebra $\cB\subset\sB(G)$ with 
$\vp_f=\vp_g\circ S$, where $f$ is cyclic for $\cA$ and $g$ is cyclic for 
$\cB'$. Then there exists a normal unital completely positive extension 
$Z\colon\sB(F)\rightarrow\sB(G)$ of $S\colon\cA\rightarrow\cB$ such that 
also the equation $\vp_f=\vp_g\circ S$ extends to $\vp_f=\vp_g\circ Z$.
\end{theorem}

\noindent
\sc{Proof.~}
We have the duality $S \leftrightarrow S'$. By Theorem \ref{mthm} there exists a weak tensor dilation for $S'$. Now the construction $j\to Z$ above or its analogue in \cite{Goh04,Goh03a} applies and yields an extension $Z$ with the required properties.~\qedsymbol

\lf
The study of extensions of completely positive maps originated in the work of Arveson \cite{Arv69}. See also Effros and Ruan \cite{EfRu00}, Chapter 4, for a recent account. Our result shows that preservation of states can be included in the extension. This is interesting, especially, when these mappings are interpreted as transition operators in noncommutative probability theory. See \cite{Goh04} for applications along these lines.


\newcommand{\Swap}[2]{#2#1}\newcommand{\Sort}[1]{}
\providecommand{\bysame}{\leavevmode\hbox to3em{\hrulefill}\thinspace}
\providecommand{\MR}{\relax\ifhmode\unskip\space\fi MR }
\providecommand{\MRhref}[2]{%
  \href{http://www.ams.org/mathscinet-getitem?mr=#1}{#2}
}
\providecommand{\href}[2]{#2}

\noindent
Rolf Gohm, \it{Institut f{\"u}r Mathematik und Informatik,
Ernst-Moritz-Arndt-Universit\"at Greifswald, Jahnstra{\ss}e 15a, 17487
Greifswald, Germany}\\E-mail: \tt{gohm@uni-greifswald.de}

\lf\noindent
Michael Skeide, \it{Dipartimento S.E.G.e S., Universit\`a degli Studi del Molise, Via de Sanctis, 86100 Campobasso, Italy}\\E-mail: \tt{skeide@math.tu-cottbus.de}

\end{document}